\documentclass{ifacconf}

\usepackage{graphicx}      
\usepackage{natbib}        

\usepackage{array}
 \usepackage{multirow}
 \usepackage{diagbox}
	\usepackage{makecell}
	\usepackage{slashbox}
	\usepackage{tabularx,array,booktabs,calc,multirow}				
	\newcolumntype{Y}{>{\centering\arraybackslash}X}	
	
\usepackage{amsmath}
\usepackage{amssymb}
\usepackage{siunitx}
\usepackage{xcolor}
\usepackage{algorithm}
\usepackage{algpseudocode}
\begin{document}
	\begin{frontmatter}
		\title{Adaptive Koopman-Based Models for Holistic Controller and Observer Design\thanksref{footnoteinfo}} 
		
		\thanks[footnoteinfo]{This research was funded by the Federal Ministry of Education and Research of Germany under the grant number 01IS20052 (DART) and the Ministry of Culture and Science of the State of North Rhine-Westphalia, Germany under the grant number NW21-059A (SAIL). The responsibility for the content of this publication lies with the authors.\\
		© 2023 the authors. \textit{This work has been accepted to IFAC for publication under a Creative Commons Licence CC-BY-NC-ND}.}
		\author{Annika Junker,}
		\author{Keno Pape,}
		\author{Julia Timmermann,}
		\author{Ansgar Trächtler}
		
		\address{Heinz Nixdorf Institute, Paderborn University, Germany (e-mail: \{annika.junker, keno.pape, julia.timmermann, ansgar.traechtler\}@hni.upb.de)}
		
		\begin{abstract}                
			We present a method to obtain a data-driven Koopman operator-based model that adapts itself during operation and can be straightforwardly used for the controller and observer design. The adaptive model is able to accurately describe different state-space regions and additionally consider unpredictable system changes that occur during operation. Furthermore, we show that this adaptive model is applicable to state-space control, which requires complete knowledge of the state vector. For changing system dynamics, the state observer therefore also needs to have the ability to adapt. To the best of our knowledge, there have been no approaches presently available that holistically use an adaptive Koopman-based plant model for the state-space design of both the controller and observer. We demonstrate our method on a test rig: controller and observer adequately adapt during operation, so that outstanding control performance is achieved even in the case of strong occuring systems changes.
		\end{abstract}
		
		
	\end{frontmatter}
	
	\section{Introduction}\label{sec:introduction}
	Data-driven methods are increasingly used for modeling technical systems due to their high performance. However, in the field of control engineering, the established design methods usually require analytical physics-derived models, so it is essential to merge these two areas. Challenging for mechatronic systems is that unpredictable changes can occur over the product life cycle, e.g., the installation of a new component, wear or temperature fluctuations during operation. This requires the information processing system to measure these changes and constantly update the internal model used for the controller and observer. 
	
	\begin{figure*}
		\begin{center}
			\includegraphics[width=1\textwidth]{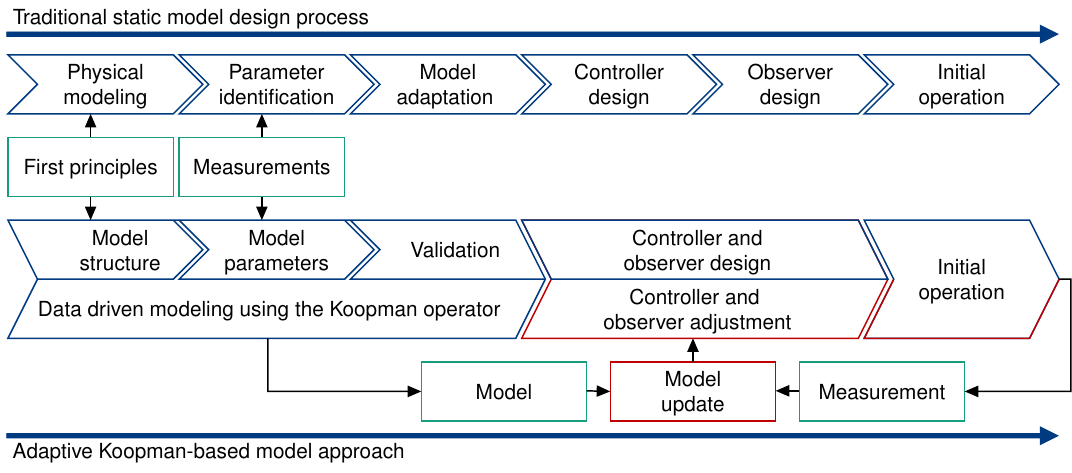}    
			\caption{Our method enables to obtain a data-driven Koopman-based model that adapts itself during operation and can be straightforwardly used for the controller and observer adjustment. The traditional approach uses the measured data only for parameter identification. In contrast, the adaptive approach requires the continuous recording and analysis of the measured data. Following the classical system identification, we exploit the promising idea of the Koopman operator, which allows to describe nonlinear dynamics in a higher-dimensional space by a linear model.}
			\label{fig:overview}
		\end{center} 
	\end{figure*}
	Research in the field of the Koopman operator provides promising strategies for model building that is appropriate from a control engineering point of view, see \cite{JTT22}, while at the same time utilizing the strengths of data-driven approaches. The underlying idea is that nonlinear dynamics is approximated by a linear but higher-dimensional operator (\cite{WKR15}), enabling the control design for a nonlinear system using linear design strategies, e.g., model predictive control (MPC) and linear-quadratic regulator (LQR) design, see \cite{KM18} and \cite{MCTM19}, respectively. The associated method is called Extended Dynamic Mode Decomposition (EDMD). 
	
	In the field of adaptive control, there are well-established methods to deal with changing system parameters, e.g., unpredictable fluctuations in operating or environmental conditions, or if an insufficiently accurate linear model is assumed to describe nonlinear dynamics. Direct adaptive control means adjusting the controller parameters in real-time, e.g., depending on the error between the system output and the output of the desired reference model, so is called \textit{MRAC} (Model Reference Adaptive Control, see \cite{Ast08}). In contrast, indirect adaptive control is a two-step process: first, the parameters of the plant model are estimated depending on the error between the plant output and the plant model output, and then the controller parameters are calculated accordingly. Since we value interpretable models with deep insights into system dynamics, e.g., the state, and want to use our adaptive internal plant model for both controller and observer design, we aim for indirect adaptive control approaches.
	
	Several approaches to adaptive Koopman-based models already exist to design the controller. \cite{PK18b}, \cite{ZRDC19} and \cite{CDGW20} utilize a recursive method to adapt a Koopman-based model online. To overcome the need for a completely measurable state, the authors in \cite{CDGW20} propose to use an input-output Koopman-based model. \cite{CSOW21} use an adaptive Koopman-based model with a velocity-based linearization to obtain a quasi-linear parameter-varying model, which is simulatively tested within a Koopman-MPC. To account for uncertainties in modeling and data, \cite{ZZL+22} develop a recursive set-membership identification, applying this method to approximate a Koopman-based model within a Lyapunov-based MPC design and simulatively showing that model adaptation improves prediction accuracy and control performance. \cite{GWJ22} use a real-time strategy for combinatorial Koopman-based modeling through iterative learning. The authors obtain a data-driven Koopman-based model for a linear and nonlinear subspace, which is subsequently used for adaptive LQR synthesis. Simulation studies show improvement in model accuracy and control performance. In contrast to the previous parameter adaptation methods, \cite{GSM21} use an improved subspace method to obtain a time-variant approximation of the Koopman operator for the highest possible prediction quality.
	
	In a state-space control scheme, all states of the system are fed back. As usually not all states can be directly measured, a state observer is required. There are only a few results published to design an observer using a Koopman-based model. The existing methods are based on static models identified offline and are furthermore solely used for open-loop prediction. \cite{SB16} and \cite{Sur16} describe a design strategy for Luenberger or Kalman-like linear observers applied to autonomous and actuated nonlinear systems. Using DMD and EDMD, respectively, \cite{ISA+15} and \cite{NM18} identify a static Koopman-based model for use within a Kalman filter to improve the prediction accuracy of the model with recent measurement data. 
		
	System changes require that not only the controller but also the observer is adaptive. To the best of our knowledge, there are no approaches presently available to deploy adaptive Koopman-based models in an observer. In this paper, we describe a method to design an adaptive Koopman-based model that can be holistically used online in a controller and observer. This control design approach thus represents an alternative to the traditional design process on the basis of a static model, where no feedback is provided during operation, see Fig.~\ref{fig:overview}. The work is divided into three main parts. First, we set the background for system identification using the Koopman operator, cf. Sec.~\ref{sec:koopman}. Then, we describe the implementation of an adaptive Koopman-based model for changing system dynamics, cf. Sec.~\ref{sec:RKM}. Finally, we demonstrate the exploitation of this model for the control design of a mechatronic system, which requires both a controller and an observer, cf. Sec.~\ref{sec:RKM_design}. We conclude with a summary and an outlook, cf. Sec.~\ref{sec:conclusion}.
	
	\section{Koopman-based System Description}\label{sec:koopman}
	Inspired by the original idea of the Koopman operator, the dynamics of nonlinear systems can be linearly approximated by \textit{lifting} the states into a higher-dimensional space, referring to \textit{Extended Dynamic Mode Decomposition (\mbox{EDMD})} (\cite{WKR15}). Below is a brief description of the algorithm; a detailed utilization of the method illustrated with examples is given in \cite{JTT22}. 
	
	In the following, continuous-time control-affine systems
	\begin{equation}
		\dot{\vec{x}}=\vec{f}(\vec{x})+\vec{B}\vec{u}
		\label{eq:control-affine_system}
	\end{equation}
	with state $\vec{x}\in\mathbb{R}^n$, input $\vec{u}\in\mathbb{R}^p$ and a constant input matrix $\vec{B}\in\mathbb{R}^{n\times p}$ are considered. For the linear Koopman-based approximation of the dynamics, $N$ observable functions $\vec{\Psi}(\vec{x})=\begin{bmatrix}
		\psi_1(\vec{x}),\psi_2(\vec{x}),\cdots,\psi_N(\vec{x})
	\end{bmatrix}^\top$ are defined, which lift the states into the higher-dimensional space. It may be useful to exploit the time derivatives of the original nonlinear dynamics for this purpose, if they are known. The algorithm approximates the dynamics of the lifted states $\vec{\Psi}(\vec{x})$ as a discrete-time system 
	\begin{equation}
		\vec{\Psi}(\vec{x}_{k+1})\approx\vec{K}_t\vec{\Psi}(\vec{x}_k)+\vec{B}_t\vec{u}_k=\begin{bmatrix}
			\vec{K}_t,\vec{B}_t
		\end{bmatrix}\begin{bmatrix}
			\vec{\Psi}(\vec{x}_k)\\ \vec{u}_k
		\end{bmatrix},
		\label{eq:edmd}
	\end{equation}
	where the index $t$ denotes the discrete-time property.
	
	With $M$ sequential measurement points
	\begin{subequations}\label{eq:data_edmdc}
		\begin{align}
			\label{eq:snapshots_x_edmdc}
			\vec{X}&=\begin{bmatrix}
				\vec{x}_1, \vec{x}_2,\cdots,\vec{x}_{M-1}
			\end{bmatrix}\in \mathbb{R}^{n\times (M-1)},\\
			\vec{X}'&=\begin{bmatrix}
				\vec{x}_2, \vec{x}_3,\cdots,\vec{x}_M
			\end{bmatrix}\in \mathbb{R}^{n\times (M-1)},\\
			\vec{U}&=\begin{bmatrix}
				\vec{u}_1,\vec{u}_2\cdots,\vec{u}_{M-1}
			\end{bmatrix}\in\mathbb{R}^{p\times(M-1)}
		\end{align}
	\end{subequations}
	for state and input results
	\begin{align}
		\vec{\Psi}(\vec{X}')&\approx\vec{K}_t\vec{\Psi(X)}+\vec{B}_t\vec{U}
		=\begin{bmatrix}
			\vec{K}_t,\vec{B}_t
		\end{bmatrix}\begin{bmatrix}
			\vec{\Psi}(\vec{X})\\ \vec{U}
		\end{bmatrix}\\
		&\Rightarrow\begin{bmatrix}
			\vec{K}_t,\vec{B}_t
		\end{bmatrix}\approx\vec{\Psi}(\vec{X}')\begin{bmatrix}
			\vec{\Psi}(\vec{X})\\\vec{U}
		\end{bmatrix}^+,
	\end{align}
	where ${\vec{K}_t\in\mathbb{R}^{N\times N}}$ is the approximated Koopman operator, ${\vec{B}_t \in\mathbb{R}^{N\times p}}$ the lifted input matrix and $^+$ denotes  the Moore-Penrose inverse matrix. 
	The resulting discrete-time system description for EDMD prediction is given by
	\begin{equation}\label{eq:edmdc_pred}
		\vec{\hat{\Psi}}(\vec{x}_{k+1})=\vec{K}_t\vec{\Psi}(\vec{x}_k)+\vec{B}_t\vec{u}_k.	
	\end{equation}
	The hat on the symbols emphasizes that the quantities are estimated, cf.~\eqref{eq:edmd}, as the observable functions usually do not span a Koopman invariant subspace.
	
	\section{Recursive Koopman Model for Varying System Dynamics}\label{sec:RKM}
	Next, we derive our algorithm for a recursive Koopman-based plant model, which we refer to as \textit{rEDMD}. The main procedure, cf. Sec.~\ref{subsec:update_rules}, has been derived by \cite{CSOW21} analogously to the classical recursive least squares method by \cite{IM11}. Due to this analogy, different modifications of the classical recursive least squares method can be integrated into the adaptation algorithm developed in this paper to ensure a stable and reasonable adaptation process, cf. Sec.~\ref{subsec:extended_algorithm}.
	\subsection{Model Update Rules}\label{subsec:update_rules}
	The goal is to update the model using continuously recorded measurement data. For this purpose, analogous to \eqref{eq:data_edmdc}, the snapshot matrices are extended in each time step 
	\begin{subequations}\label{eq:data_extended}
		\begin{align}
			\label{eq:snapshots_x_extended}
			\vec{X}_k&=\begin{bmatrix}
				\vec{x}_1,\vec{x}_2,\cdots,\vec{x}_{k-1}
			\end{bmatrix}\in \mathbb{R}^{n\times (k-1)},\\
			\vec{X}'_k&=\begin{bmatrix}
				\vec{x}_2,\vec{x}_3,\cdots,\vec{x}_k
			\end{bmatrix}\in \mathbb{R}^{n\times (k-1)},\\
			\vec{U}_k&=\begin{bmatrix}
				\vec{u}_1,\vec{u}_2,\cdots,\vec{u}_{k-1}
			\end{bmatrix}\in \mathbb{R}^{p\times (k-1)},\\
			\vec{X}_{k+1}&=\begin{bmatrix}
				\vec{X}_k, \vec{x}_k
			\end{bmatrix}\in \mathbb{R}^{n\times k},\\
			\vec{X}'_{k+1}&=\begin{bmatrix}
				\vec{X}'_k, \vec{x}_{k+1}
			\end{bmatrix}\in \mathbb{R}^{n\times k},\\
			\vec{U}_{k+1}&=\begin{bmatrix}
				\vec{U}_k,\vec{u}_k
			\end{bmatrix}\in\mathbb{R}^{p\times k},
		\end{align}
	\end{subequations}
	enabling the Koopman-based model to be calculated
	\begin{align}
		\begin{split}\label{eq:redmd_k}
			\begin{bmatrix}
				\vec{K}_{t,k},\vec{B}_{t,k}
			\end{bmatrix}&=\vec{\Psi}(\vec{X}'_k)\begin{bmatrix}
				\vec{\Psi}(\vec{X}_k)\\\vec{U}_k
			\end{bmatrix}^+
		\end{split},\\
		\begin{split}\label{eq:redmd_k+1}
			\begin{bmatrix}
				\vec{K}_{t,k+1},\vec{B}_{t,k+1}
			\end{bmatrix}&=\vec{\Psi}(\vec{X}'_{k+1})\begin{bmatrix}
				\vec{\Psi}(\vec{X}_{k+1})\\\vec{U}_{k+1}
			\end{bmatrix}^+,
		\end{split}
	\end{align}
	where $k$ is the current time step. With \eqref{eq:moore-penrose} follows
	\begin{align}
		\begin{split}\label{eq:redmd_k_Gamma}
			\begin{bmatrix}
				\vec{K}_{t,k},\vec{B}_{t,k}
			\end{bmatrix}&=\vec{\Psi}(\vec{X}'_{k})\begin{bmatrix}
				\vec{\Psi}^\top(\vec{X}_{k}),\vec{U}^\top_{k}
			\end{bmatrix}\vec{\Gamma}_{k},
		\end{split}\\
		\begin{split}\label{eq:redmd_k+1_Gamma}
			\begin{bmatrix}
				\vec{K}_{t,k+1},\vec{B}_{t,k+1}
			\end{bmatrix}&=\vec{\Psi}(\vec{X}'_{k+1})\begin{bmatrix}
				\vec{\Psi}^\top(\vec{X}_{k+1}),\vec{U}^\top_{k+1}
			\end{bmatrix}\vec{\Gamma}_{k+1}
		\end{split}
	\end{align}
	with
	\begin{align}
		\begin{split}\label{eq:gamma_k}
			\vec{\Gamma}_{k}=&\left(\begin{bmatrix}
				\vec{\Psi}(\vec{X}_k)\\\vec{U}_k
			\end{bmatrix}\begin{bmatrix}
				\vec{\Psi}^\top(\vec{X}_k),\vec{U}^\top_k\\
			\end{bmatrix}\right)^{-1}\\
			=&\begin{bmatrix}
				\vec{\Psi}(\vec{X}_k)\vec{\Psi}^\top(\vec{X}_k)&
				\vec{\Psi}(\vec{X}_k)\vec{U}^\top_k\\
				\vec{U}_k\vec{\Psi}^\top(\vec{X}_k)&
				\vec{U}_k\vec{U}^\top_k
			\end{bmatrix}^{-1},
		\end{split}\\
		\begin{split}\label{eq:gamma_k+1}
			\vec{\Gamma}_{k+1}=&\left(\begin{bmatrix}
				\vec{\Psi}(\vec{X}_k)&\vec{\Psi}(\vec{x}_k)\\
				\vec{U}_k&\vec{u}_k
			\end{bmatrix}\begin{bmatrix}
				\vec{\Psi}^\top(\vec{X}_k)&\vec{U}^\top_k\\
				\vec{\Psi}^\top(\vec{x}_k)&\vec{u}^\top_k
			\end{bmatrix}\right)^{-1}\\
			=&\left[\begin{matrix}
				\vec{\Psi}(\vec{X}_k)\vec{\Psi}^\top(\vec{X}_k)+\vec{\Psi}(\vec{x}_k)\vec{\Psi}^\top(\vec{x}_k)\\
				\vec{U}_k\vec{\Psi}^\top(\vec{X}_k)+\vec{u}_k\vec{\Psi}^\top(\vec{x}_k)
			\end{matrix}\right.\cdots\\
			&\cdots\left.\begin{matrix}
				\vec{\Psi}(\vec{X}_k)\vec{U}^\top_k+\vec{\Psi}(\vec{x}_k)\vec{u}^\top_k\\
				\vec{U}_k\vec{U}^\top_k+\vec{u}_k\vec{u}^\top_k
			\end{matrix}\right]^{-1},
		\end{split}
	\end{align}
	where $\vec{\Gamma}_k,\vec{\Gamma}_{k+1}\in\mathbb{R}^{(N+p)\times(N+p)}$ are the covariance matrices for different time steps $k$, $k+1$. Combining \eqref{eq:gamma_k}-\eqref{eq:gamma_k+1} yields
	\begin{equation}\label{eq:gammas}
		\vec{\Gamma}_{k+1}^{-1}=\vec{\Gamma}_k^{-1}+\begin{bmatrix}
			\vec{\Psi}(\vec{x}_k)\\\vec{u}_k
		\end{bmatrix}\begin{bmatrix}
			\vec{\Psi}^\top(\vec{x}_k),\vec{u}^\top_k
		\end{bmatrix}
	\end{equation}
	and with \eqref{eq:redmd_k_Gamma} 
	\begin{equation}
		\begin{bmatrix}\label{eq:KBK_Gamma_k_-1}
			\vec{K}_{t,k},\vec{B}_{t,k}
		\end{bmatrix}\vec{\Gamma}_k^{-1}
		=\vec{\Psi}(\vec{X}'_{k})\begin{bmatrix}
			\vec{\Psi}^\top(\vec{X}_{k}),\vec{U}^\top_{k}
		\end{bmatrix}.
	\end{equation}
	Reshaping \eqref{eq:redmd_k+1_Gamma} yields
	\begin{align}
		\begin{split}
			\begin{bmatrix}
				\vec{K}_{t,k+1},\vec{B}_{t,k+1}
			\end{bmatrix}=&\vec{\Psi}(\vec{X}'_k)\begin{bmatrix}
				\vec{\Psi}^\top(\vec{X}_k),\vec{U}^\top_k
			\end{bmatrix}\vec{\Gamma}_{k+1}\\
			&+\vec{\Psi}(\vec{x}_{k+1})\begin{bmatrix}
				\vec{\Psi}^\top(\vec{x}_k),\vec{u}^{\top}_k
			\end{bmatrix}\vec{\Gamma}_{k+1},
		\end{split}
	\end{align}
	and with \eqref{eq:gammas}-\eqref{eq:KBK_Gamma_k_-1} follows
	\begin{align}
	\begin{split}
		&\begin{bmatrix}
			\vec{K}_{t,k+1},\vec{B}_{t,k+1}
		\end{bmatrix}=\begin{bmatrix}
			\vec{K}_{t,k},\vec{B}_{t,k}
		\end{bmatrix}\\
		&+\left(\vec{\Psi}(\vec{x}_{k+1})-\begin{bmatrix}\vec{K}_{t,k},\vec{B}_{t,k}
		\end{bmatrix}\begin{bmatrix}
			\vec{\Psi}(\vec{x}_k)\\\vec{u}_k
		\end{bmatrix}\right)\vec{\gamma}_k,
	\end{split}
\end{align}
	where 
	\begin{equation}\label{eq:correction_vector_short}
		\vec{\gamma}_k=\begin{bmatrix}\vec{\Psi}^\top(\vec{x}_k),\vec{u}_k^\top\end{bmatrix}\vec{\Gamma}_{k+1}
	\end{equation}
	is the correction vector and $\vec{\Gamma}_{k+1}$ is calculated from \eqref{eq:gammas} with \eqref{eq:matrix_inversion}, yielding
	\begin{equation}\label{eq:gamma_k+1_recursive_long}
		\begin{split}
			\vec{\Gamma}_{k+1}=\vec{\Gamma}_{k}-\frac{\vec{\Gamma}_{k}\begin{bmatrix}
					\vec{\Psi}(\vec{x}_k) \\ \vec{u}_k
				\end{bmatrix}\begin{bmatrix}
					\vec{\Psi}^\top(\vec{x}_k), \vec{u}_{k}^\top
				\end{bmatrix}\vec{\Gamma}_{k}}{\begin{bmatrix}
					\vec{\Psi}^{\top}(\vec{x}_k),\vec{u}^{\top}_k
				\end{bmatrix}\vec{\Gamma}_{k}\begin{bmatrix}
					\vec{\Psi}(\vec{x}_k) \\ \vec{u}_k
				\end{bmatrix} +1}
		\end{split}
	\end{equation}
	and with \eqref{eq:gamma_k+1_recursive_long} results
	\begin{equation}\label{eq:correction_vector_long}
			\vec{\gamma}_{k}= \frac{\begin{bmatrix}
					\vec{\Psi}^{\top}(\vec{x}_k),\vec{u}^{\top}_k
				\end{bmatrix}\vec{\Gamma}_{k}}{\begin{bmatrix}
					\vec{\Psi}^{\top}(\vec{x}_k),\vec{u}^{\top}_k
				\end{bmatrix}\vec{\Gamma}_{k}\begin{bmatrix}
					\vec{\Psi}(\vec{x}_k) \\ \vec{u}_k
				\end{bmatrix} +1}
	\end{equation}
	and 
	\begin{equation}\label{eq:gamma_k+1_recursive_short}
		\vec{\Gamma}_{k+1}=\vec{\Gamma}_{k} \left(\vec{I} - \begin{bmatrix}
			\vec{\Psi}^{\top}(\vec{x}_k),\vec{u}^{\top}_k
		\end{bmatrix}\vec{\gamma}_k\right).
	\end{equation}
	
	We now introduce the forgetting factor $\lambda$, which enables past measurements to be \textit{forgotten} in the event of system changes, where the smaller $\lambda$ is, the less influence past measurements have on the current model estimate. Thus, the adaptation speed can be systematically controlled. $\lambda$ affects the calculation of $\gamma_k$ and $\Gamma_{k+1}$, resulting in the equations for the model update
	\begin{align}
		\vec{\gamma}_{k}= \frac{\begin{bmatrix}
				\vec{\Psi}^{\top}(\vec{x}_k),\vec{u}^{\top}_k
			\end{bmatrix}\vec{\Gamma}_{k}}{\begin{bmatrix}
				\vec{\Psi}^{\top}(\vec{x}_k),\vec{u}^{\top}_k
			\end{bmatrix}\vec{\Gamma}_{k}\begin{bmatrix}
				\vec{\Psi}(\vec{x}_k) \\ \vec{u}_k
			\end{bmatrix} +\lambda},\label{eq:alg1}
	\end{align}
	\begin{align}
		\begin{split}\label{eq:alg2}
			&\begin{bmatrix}
				\vec{K}_{t,k+1},\vec{B}_{t,k+1}
			\end{bmatrix}=\begin{bmatrix}
				\vec{K}_{t,k},\vec{B}_{t,k}
			\end{bmatrix}\\
			&+\left(\vec{\Psi}(\vec{x}_{k+1})-\begin{bmatrix}\vec{K}_{t,k},\vec{B}_{t,k}
			\end{bmatrix}\begin{bmatrix}
				\vec{\Psi}(\vec{x}_k)\\\vec{u}_k
			\end{bmatrix}\right)\vec{\gamma}_k,
		\end{split}\\
		&\vec{\Gamma}_{k+1}=\vec{\Gamma}_{k} \left(\vec{I} - \begin{bmatrix}
			\vec{\Psi}^{\top}(\vec{x}_k),\vec{u}^{\top}_k
		\end{bmatrix}\vec{\gamma}_k\right)\frac{1}{\lambda}\label{eq:alg3}.
	\end{align}
	Thus, \eqref{eq:alg1}, \eqref{eq:alg2} and \eqref{eq:alg3} form the core of our algorithm.

	\subsection{Extended Algorithm}\label{subsec:extended_algorithm}
	Our algorithm includes certain extensions to ensure that the adaption process is stable and reasonable. Below, we present the different extensions resulting in our final algorithm.
	\subsubsection{Checking the Need to Update the Model}\label{subsubsec:update}
	To decide whether the model needs to be updated, the prediction accuracy of the current model serves as a criterion, as proposed by \cite{CDGW20}. During operation, $M_{op}$ measurements of the states and the inputs are stored
	\begin{subequations}\label{eq:data_edmdc_m_op}
	\begin{align}
		\bar{\vec{X}}_{k}^{\top} &= \begin{bmatrix} \vec{x}_{k-M_{op}}^{\top}, \vec{x}_{k-M_{op}+1}^{\top}, \cdots, \vec{x}_{k}^{\top} \end{bmatrix},\\
		\bar{\vec{U}}_{k}^{\top} &= \begin{bmatrix} \vec{u}_{k-M_{op}}^{\top}, \vec{u}_{k-M_{op}+1}^{\top}, \cdots, \vec{u}_{k}^{\top} \end{bmatrix}.
	\end{align}
	\end{subequations}
	The current EDMD model is simulated according to \eqref{eq:edmdc_pred} and the states $\hat{\vec{x}}$ are extracted from the lifted states with a projection matrix $\vec{P}_x$
	\begin{align}
		\hat{\vec{x}}_{k-M_{op}+1} &= \vec{P}_x\left(\vec{K}_{t,k}\vec{\Psi}(\vec{x}_{k-M_{op}}) + \vec{B}_{t,k}\vec{u}_{k-M_{op}}\right),\\
		\hat{\vec{x}}_{k-M_{op}+2} &= \vec{P}_x\left(\vec{K}_{t,k}\vec{\Psi}(\hat{\vec{x}}_{k-M_{op}+1}) + \vec{B}_{t,k}\vec{u}_{k-M_{op}+1}\right),\\
		&\hspace{0.2cm}\vdots\notag\\
		\hat{\vec{x}}_{k} &= \vec{P}_x\left(\vec{K}_{t,k}\vec{\Psi}(\hat{\vec{x}}_{k-1}) + \vec{B}_{t,k}\vec{u}_{k-1}\right),
\end{align}
	resulting in
	\begin{equation}
		\hat{\vec{X}}_{k}^{\top} = \begin{bmatrix} \vec{x}_{k-M_{op}}^{\top}, 
		\hat{\vec{x}}_{k-M_{op}+1}^{\top} ,\cdots, \hat{\vec{x}}_{k}^{\top} \end{bmatrix}.
	\end{equation}
	
	If the defined accuracy limit $\varepsilon_{low}$ is exceeded, i.e.,
	\begin{equation}\label{eq:acc_limit}
		\left\Vert \bar{\vec{X}}_k - \hat{\vec{X}}_k \right\Vert_{\infty} \geq \varepsilon_ {low},
	\end{equation}
	the new measurement $\begin{bmatrix}
		\vec{\Psi}^\top(\vec{x}_{k+1}), \vec{u}_{k+1}^\top
	\end{bmatrix}$ is used for a model update.  
	
	\subsubsection{Varying Forgetting Factor}\label{subsubsec:forgetting_factor}
	Due to the better performance, we use a variable forgetting factor $\lambda_k$, see \cite{FKY81},
	\begin{equation}\label{eq:var_lambda}
		\begin{split}
			\lambda_{k+1} &= 1 - \frac{1}{\Sigma_{0,k}}\left( 1 - \begin{bmatrix}
				\vec{\Psi}^{\top}(\vec{x}_k), \vec{u}_k^{\top}
			\end{bmatrix}\vec{\gamma}_k^{\top} \right) e^2_{post,k},\\ 
		\end{split}
	\end{equation}
	where 
	\begin{equation}\label{eq:e_post}
		e_{post,k}= y_k - \vec{P}_y\left(\vec{K}_{t,k}\vec{\Psi}(\vec{x}_{k-1}) + \vec{B}_{t,k}\vec{u}_{k-1}\right)
	\end{equation}
	is the error between the measured and the predicted system output. $\vec{P}_y$ is a projection matrix to extract the output $y$ from the lifted states and $\Sigma_{0,k}=\sigma_{e,k} N_0$ is an assumed measurement noise with the variance $\sigma_{e,k}$ , which is a posteriori determined by the error behavior \eqref{eq:e_post}, and $N_0$ is an adjustable sensitivity factor. If the accuracy, cf.~\eqref{eq:acc_limit}, even exceeds $\varepsilon_{high}$, the adaptation speed is increased by multiplying $\Sigma_{0,k}$ by a gain factor $ \mu_{\Sigma}>1$. Note here that the described strategy is only applicable for systems with a one-dimensional output, otherwise the use of multiple varying forgetting factors would be necessary. 
	
	\subsubsection{Setting the Trace}\label{subsubsec:constant_trace}
	Using a forgetting factor can lead to a covariance windup in operating situations with insufficient excitation because the elements of the covariance matrix $\vec{\Gamma}_k$ may strongly increase, causing the estimation algorithm to react sensitively to changes in the measured data. Therefore, following \cite{LG85}, we exploit a constant-trace algorithm for the covariance matrix to limit the sensitivity and thus the adaptation intensity of our algorithm. In case the current trace of the covariance matrix exceeds a maximum value, the elements of the covariance matrix are multiplied by a reduction factor $\mu_{\vec{\Gamma}}<1$
	\begin{equation} \label{eq:const_trace_algo}
		\begin{split}
			\vec{\Gamma}_{k} \leftarrow \vec{\Gamma}_{k}\mu_{\vec{\Gamma}} \text{  with  }
			\mu_{\vec{\Gamma}} = \frac{\text{tr}\left(\vec{\Gamma}_{max}\right)}{\text{tr}\left(\vec{\Gamma}_{k}\right)}.
		\end{split}
	\end{equation} 
	
	Algorithm~\ref{alg:redmd} summarizes rEDMD with extensions.
	\algnewcommand{\Initialize}[1]{%
		\State \textbf{Initialize: }\parbox[t]{.6\linewidth}{\raggedright #1}
	}
	\begin{algorithm}[!ht]
		\begin{algorithmic}
			\Initialize{$\vec{K}_{t,0},\vec{B}_{t,0},\vec{\Gamma}_{0}, \lambda_0, \Sigma_{0,0}$}
			\While{$k \leq T_{sim}$} 
			\If{$M_{op}$ samples of data were collected}
			\If{$\left\Vert \bar{\vec{X}}_k - \hat{\vec{X}}_k \right\Vert_{\infty} \geq \varepsilon_ {low}$}
			\If{$\text{tr}(\vec{\Gamma}_k) > \text{tr}(\vec{\Gamma}_{max})$}
			\State Set $\vec{\Gamma}_k$ with \eqref{eq:const_trace_algo}
			\EndIf
			\If{$\left\Vert \bar{\vec{X}}_k - \hat{\vec{X}}_k \right\Vert_{\infty} \geq \varepsilon_ {high}$}					
			\State Set $\Sigma_{0,k} \leftarrow \Sigma_{0,k} \mu_{\Sigma}$
			\EndIf
			\State{Update rEDMD using \eqref{eq:alg1}, \eqref{eq:alg2} and \eqref{eq:alg3}}
			\State{Calculate $\lambda_{k+1}$ with \eqref{eq:var_lambda}}
			\Else
			\State{No Update of rEDMD and $\lambda$}
			\EndIf
			\Else
			\State{Update rEDMD with inital values and \eqref{eq:alg1}-\eqref{eq:alg3}}
			\EndIf				
			\EndWhile
		\end{algorithmic}
		\caption{rEDMD algorithm with extensions}\label{alg:redmd}
	\end{algorithm}

	\section{Results for a Holistic Adaptive Controller and Observer Design}\label{sec:RKM_design}
	For changing systems where not all states are known and state feedback control is used, it is essential to design both the controller and the required observer as adaptive systems. This requires a holistic strategy, i.e., both designs must be equally regarded. The linearity of the rEDMD model makes it possible to systematically employ existing linear design procedures. Thus, we propose to use a model predictive control in conjunction with a Kalman filter. Note here that the linear design process is realized in a higher-dimensional space, i.e., with more states than the original system, so is $N>n$. 
	
	We experimentally demonstrate the success of our holistic method on the golf robot, whose goal is to autonomously learn to hit a golf ball into the hole from an arbitrary starting point on the green using combined data-driven and physics-based methods. A detailed description of the golf robot is given in~\cite{JFTT22b}, where the overall problem is decomposed into seperate sub-tasks. The stroke speed control of the golf club is a challenging task in terms of control engineering because of nonlinear effects, changes in dynamics due to the replacement of the golf club and additionally because not all states can be directly measured. Based on these considerations, we selected this subtask to demonstrate the success of our method. 
	
	\begin{figure}[ht!]
		\begin{center}
			\includegraphics[width=0.45\textwidth]{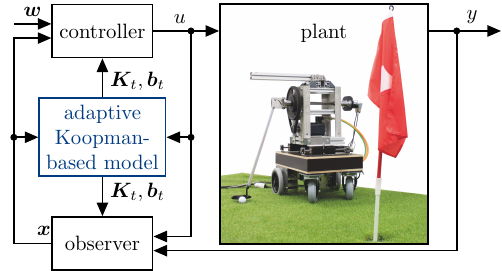}    
			\caption{We show how the adaptive Koopman-based model is used for both the controller and the observer design considering the stroke mechanism of our golf robot. The lower gear shaft is driven with input torque $u$, while the golf club is mounted on the upper gear shaft with the state vector  ${\vec{x}=\left[x_1,x_2\right]^\top}$ containing the angular position and angular velocity. $\vec{w}$ specifies the reference trajectories for $\vec{x}$ and the measurable system output is $y=x_1$.}
			\label{fig:control_loop}
		\end{center} 
	\end{figure}
\begin{figure*}
		\begin{center}
			\includegraphics[width=1\textwidth]{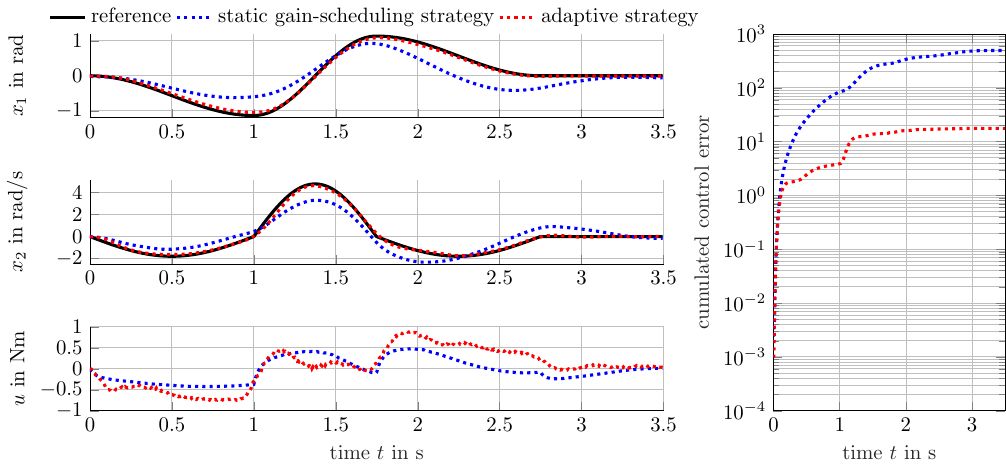}    
			\caption{The control task of the stroke mechanism is to precisely hit the ball with a desired stroke velocity. For this purpose, three phases \textit{lunge, strike, reset} are defined so that the club has its maximum desired angular velocity $x_2$ at the angular position $x_1=0$, where the ball is supposed to be hit. The mathematical description of the reference trajectory is derived in \cite{JFTT22b}. Our strategy, which uses the adaptive rEDMD model for both the controller and the observer, provides significantly higher control performance for a system with changes than a strategy using a static gain scheduling design approach for the controller and observer. Note here that the desired stroke speed of \SI{3}{\meter\per\second} at the club corresponds to a desired angular velocity of \SI{5}{\radian\per\second} due to the club length of \SI{0.6}{\meter}.}
			\label{fig:adaptive_plot}
		\end{center}
	\end{figure*}

	The procedure for approximating a static EDMD model is described in \cite{JTT22} and extended by Algorithm~\ref{alg:redmd} to include adaptability in this work. The dimension of the original system is $n=2$ and is increased to $N=4$ by the Koopman approach. We initialized the rEDMD model offline using historical open-loop training data and implemented a model predictive control with a Kalman filter using our adaptive Koopman-based model, see Fig.~\ref{fig:control_loop}. Then we changed the robot by mounting an extra mass and thus doubling the golf club mass  and additionally by increasing the frictional torque of the driven shaft by tightening belts. Note here that these unrealistically strong system changes are for illustration purposes only to demonstrate the power of our algorithm.  Fig.~\ref{fig:adaptive_plot} illustrates the resulting control performance on the basis of the cumulated quadratic control error
	\begin{equation}\label{eq:control_error}
		e(t_k)=\sum_{m=1}^{k} \left\Vert\vec{w}(t_m)-\vec{x}(t_m)\right\Vert^2_2,
	\end{equation}
	where $\vec{w}(t)$ is the reference. In comparison, the control performance of the state-of-the-art controller, which is based on a traditional static design method with a gain-scheduling approach, cf.~\cite{JFTT22b}, is shown. Fig.~\ref{fig:adaptive_plot} illustrates an example stroke with a desired speed of \SI{3}{\meter\per\second} at the club head. However, the success of our method can also be transferred to other stroke speeds, see Table~\ref{tab:error}, and also works for the case where no system changes occur. Measurement noise influence on the identification could be reduced for the presented experiment by careful selection of the algorithm parameters. Since noise is generally a major challenge and always present, it is currently under further investigation.
	\begin{table}[h!]
		\label{tab:error}
		\caption{Normalized control errors \eqref{eq:control_error} for different stroke speeds at the golf robot without and with occuring system changes. The errors are multiplied by $10^2$ for better readability.}
		\begin{center} 
			\begin{tabularx}{\linewidth}{
					>{\setlength\hsize{0.5\hsize}}l
					>{\setlength\hsize{0.1\hsize}}X
					>{\setlength\hsize{0.1\hsize}}X
					>{\setlength\hsize{0.1\hsize}}X
					>{\setlength\hsize{0.1\hsize}}X
				} 
				& \multicolumn{2}{c}{\makecell{w/o system\\ changes}}& \multicolumn{2}{c}{\makecell{w/ system \\ changes }} \\
				\cmidrule[\heavyrulewidth](lr){1-5}
				\multicolumn{1}{l}{stroke speed $v$}& $\SI{1}{\meter\per\second}$&$\SI{3}{\meter\per\second}$ &$\SI{1}{\meter\per\second}$&$\SI{3}{\meter\per\second}$   \\ 
				\cmidrule[\heavyrulewidth](lr){1-5}\cmidrule[\heavyrulewidth](lr){2-5}
				gain-scheduling strategy & 5.68 & 6.33 & 46.00 & 55.08 \\
				\cmidrule(lr){1-5}
				\makecell[l]{adaptive controller and\\ non-adaptive observer} & 0.31 & 0.96 & 6.36 & 5.08  \\			
				\cmidrule(lr){1-5}
				\textbf{\makecell[l]{both adaptive\\controller and observer}} &\textbf{0.34} & \textbf{0.67}& \textbf{1.38} &\textbf{2.18}   \\
				\cmidrule[\heavyrulewidth](lr){1-5}%
			\end{tabularx}
		\end{center}
	\end{table}
	\section{Conclusion \& Outlook}\label{sec:conclusion}
	We established a method to design an adaptive Koopman-based model and demonstrated the success when using it as an internal model for the controller and observer in a state-space control scheme achieving outstanding control performance on a test rig. The complete observability of the specific rEDMD model for the golf robot has been verified offline. However, the generalized influence of the higher-dimensional model on the observability remains an open research question. Furthermore, it is a matter of interest how to extend our approach to structural system changes such as nonlinearities that have not occured before. For system descriptions which are valuable from a control engineering point of view, the proposed recursive parameter estimation method may be transferred to data-driven PCHD models (\cite{JTT22b}).

	\appendix
	\section{Matrix Calculation}\label{app:moore-penrose}
	The Moore-Penrose inverse matrix of a matrix $\vec{A}$ with full row rank can be calculated by
	\begin{equation}\label{eq:moore-penrose}
		\vec{A}^+=\vec{A}^\top\left(\vec{A}\vec{A}^\top\right)^{-1}.
	\end{equation}
	If $\vec{A}$, $\vec{C}$ and $\left(\vec{A}^{-1}+\vec{B}\vec{C}^{-1}\vec{D}\right)$ are non-singular quadratic matrices and 
	\begin{equation}
		\vec{E}=\left(\vec{A}^{-1}+\vec{B}\vec{C}^{-1}\vec{D}\right)^{-1}
	\end{equation}
	then
	\begin{equation}\label{eq:matrix_inversion}
		\vec{E}=\vec{A}-\vec{A}\vec{B}\left(\vec{D}\vec{A}\vec{B}+\vec{C}\right)^{-1}\vec{D}\vec{A},
	\end{equation}
	see \cite{IM11}.
	
\end{document}